\documentclass[psamsfonts,reqno]{amsart}

\usepackage{amssymb}
\usepackage{amscd}
\usepackage{eucal}
\usepackage{psfig}
\usepackage[dvips]{graphicx}
\usepackage[small]{caption}

\DeclareMathOperator{\lm}{\text{l.i.m.}}

\newcommand{\beq}{\begin{equation}}
\newcommand{\eeq}{\end{equation}}
\newcommand{\sds}{\strut\displaystyle}

\newcommand{\Real}{\mbox{Re\,}}
\newcommand{\Imag}{\mbox{Im\,}}
\newcommand{\R}{\mathbb{R}}
\newcommand{\Z}{\mathbb{Z}}
\newcommand{\C}{\mathbb{C}}
\newcommand{\SSS}{\mathbb{S}}

\newcommand{\cA}{\mathcal{A}}
\newcommand{\cB}{\mathcal{B}}
\newcommand{\cC}{\mathcal{C}}
\newcommand{\cF}{\mathcal{F}}
\newcommand{\cI}{\mathcal{I}}
\newcommand{\cL}{\mathcal{L}}
\newcommand{\cM}{\mathcal{M}}

\newtheorem{theorem}{Theorem}
\newtheorem*{theorem-nonumber}{Theorem}
\newtheorem{proposition}{Proposition}

\theoremstyle{remark}
\newtheorem*{remark}{Remark}
\newtheorem*{remarks}{Remarks}

\theoremstyle{definition}

\newcommand{\sfrac}[2]{{\vphantom1\smash{\lower.5ex\hbox{\small$#1$}}\over
        \vphantom1\smash{\raise.4ex\hbox{\small$#2$}}}} 

\begin{document}

\title[Hausdorff Moments and Power Series]{Hausdorff Moments, Hardy Spaces and Power Series}

\author[E. De Micheli]{E. ~De Micheli}
\address[E. De Micheli]{IBF -- Consiglio Nazionale delle Ricerche \\ Via De Marini, 6 - 16149 Genova, Italy}
\email[E.~De Micheli]{demicheli@ge.cnr.it}

\author[G. A. Viano]{G. A. ~Viano}
\address[G. A. ~Viano]{Dipartimento di Fisica - Universit\`a di Genova,
Istituto Nazionale di Fisica Nucleare - sez. di Genova,\\
Via Dodecaneso, 33 - 16146 Genova, Italy}
\email[G.A.~Viano]{viano@ge.infn.it}

\begin{abstract}
In this paper we consider power and trigonometric series whose coefficients
are supposed to satisfy the Hausdorff conditions, which play a relevant role
in the moment problem theory. We prove that these series converge to functions
analytic in cut domains. We are then able to reconstruct the jump functions
across the cuts from the coefficients of the series expansions by the use of
the Pollaczek polynomials. We can thus furnish a solution for a class of
Cauchy integral equations.
\end{abstract}

\maketitle

\section{Introduction}
\label{introduction_section}
The problem of characterizing the analytic properties of the functions,
in terms of the coefficients of their power expansions, is very old
and goes back to a classical result due to Le Roy \cite{Leroy}:

\begin{theorem-nonumber}[Le Roy]
\label{the:Leroy}
If in the Taylor series $\sum_{n=0}^\infty a_n\,z^n$
the coefficients $a_n$ are the restriction to the integers of a function
$\tilde{a} (\lambda) (\lambda \in \C)$, holomorphic in the half--plane
$\Real \lambda > -1/2$, and moreover there exist two constants $A$ and $N$
such that
\beq
\label{zerouno}
|\tilde{a} (\lambda) | \leq A (1+|\lambda|)^N \qquad (\text{\rm\Real} \lambda > -\sfrac{1}{2}),
\eeq
then the series converges to a function $f(z)$ analytic in the unit disk $D=\{z \mid |z| < 1\}$,
and furthermore $f(z)$ admits a holomorphic extension to the cut plane
$\{z \in \C \,\mid \,[1,+\infty) \,\}$.
\end{theorem-nonumber}

Other similar results are due to Lindel\"{o}f \cite{Lindelof} and
Bieberbach \cite{Bieberbach}.
More recently, Stein and Wainger \cite{Stein} have reconsidered the problem
in the framework of the Hardy space theory. More precisely, they assume
that the coefficients $a_n$ are the restriction to the integers of a function
$\tilde{a}(\lambda)$ holomorphic in the half--plane $\Real \lambda > -1/2$, and,
in addition,
$\tilde{a}(\lambda)$ is supposed to belong to the Hardy space
$H^2(\C_{-1/2})$ with norm
$\|\tilde{a}\|_2 = \sup_{\sigma > -\frac{1}{2}} \left ( \int_{-\infty}^{\infty}
| \tilde{a} (\sigma+i\nu)|^2\, d\nu\right )^{1/2}$ $(\C_{-1/2} = \{\lambda\in\C \, |\, \Real\lambda > -1/2\})$.
Correspondingly, they consider the class of functions $f(z)$ analytic in the
complex $z$--plane slit along the positive real axis from $1$ to $+\infty$
(this domain is denoted by $\SSS$). The space of functions analytic in $\SSS$ for which
$\sup_{y\neq 0}\left ( \int_{-\infty}^{+\infty} |f(x+iy)|^2\, dx\right )^{1/2} < \infty$
is denoted by $H^2(\SSS)$. These authors have proved the following result:

\begin{theorem-nonumber}[Stein--Wainger]
\label{the:Stein-Wainger}
Suppose that $f(z)=\sum_{n=0}^\infty a_n z^n$. Then $f \in H^2(\SSS)$ if and only if
$a_n = [\tilde{a}(\lambda)]_{(\lambda=n)}$, where $\tilde{a} \in H^2(\C_{-1/2})$.
Moreover,
\beq
\label{zerodue}
\| F \|_2^2 = 2\pi \| \tilde{a}\|_2^2,
\eeq
where
\beq
\label{zerotre}
\| \tilde{a}\|_2^2 = \int_{-\infty}^{+\infty} \left|\tilde{a}(-\frac{1}{2}+i\nu)\right|^2\,d\nu,
\eeq
and
\beq
\label{zeroquattro}
\| F \|_2^2 = \int_{1}^{+\infty} |F(x)|^2\,dx,
\eeq
$F(x)$ being the jump function, i.e., $F(x)= -i(f_+(x) - f_-(x))$;
$f_+(x)$ and $f_-(x)$ are the boundary values of
$f(x+iy)$ and $f(x-iy)$, respectively
$(f_\pm (x) = \lim_{\substack{y\rightarrow 0 \\ y>0}} f(x\pm iy))$.
\end{theorem-nonumber}

The main purpose of the present paper is the reconstruction of the jump function
across the cut from the coefficients of the power expansion. This result
is achieved by the use of the Pollaczek polynomials \cite{Bateman,Szego}, and it
will be illustrated in Section \ref{sezione3}.
This reconstruction allows us to solve the Cauchy integral equation of the
following type:
\beq
\label{zeroquattroprimo}
f(z) = \frac{1}{2\pi} \int_{1}^{+\infty} \frac{F(x)}{x-z}\,dx,
\eeq
when the Taylor coefficients $a_n=f^{(n)}(0)/n!$, $(n=0,1,2,\ldots)$
are supposed to be known.
Unfortunately, in the numerical analysis and in the applications to physical problems
only a finite number of Taylor coefficients are known, and, moreover, they are affected
by noise or, at least, by round--off errors. Furthermore, the integral equations of first
kind, like equations (\ref{zeroquattroprimo}), give rise to the so--called ill--posed
problems in the sense of Hadamard \cite{Hadamard}: The solution does not depend continuously
on the data. We shall briefly return on this important point in Section \ref{sezione3}.
In a separate paper we shall discuss in detail how to manage numerically the
method presented here.

In the theorems of Leroy and Stein--Wainger the coefficients
$a_n$ are required to be the restriction of a function $\tilde{a}(\lambda) (\lambda\in\C)$
holomorphic in the half--plane $\Real\lambda > -1/2$,
and, in the case of the Stein--Wainger theorem, this function is also assumed
to belong to the Hardy space $H^2(\C_{-1/2})$.
We prefer to start by requiring that the coefficients
$a_n$ satisfy the so--called \emph{Hausdorff conditions} \cite{Widder}, which guarantee that
they can be regarded as Hausdorff moments in a sense that will be
explained in Section \ref{sezione1}. This approach is convenient for several reasons:
\begin{itemize}
\item[i)] From the Hausdorff conditions and by the use of the Carlson theorem \cite{Boas}
we derive immediately the existence of a unique interpolation of the coefficients
which is holomorphic in a half--plane.
\item[ii)] Following Watanabe \cite{Watanabe}, we can give a suggestive probabilistic
interpretation of the Hausdorff conditions. Accordingly, the coefficients $a_n$
can be regarded as the values of an harmonic function associated with a
Markov process in a sense that will be outlined in Section \ref{sezione1}.
\item[iii)] By imposing to the coefficients $a_n$ Hausdorff conditions of
various types, we can, correspondingly, obtain more specific properties of
smoothness of the function which gives the jump across the cut.
\end{itemize}
The last point is particularly relevant. In fact, in the Stein--Wainger approach
one works essentially with the unitary equivalence between $H^2(\SSS)$
and $L^2(1,+\infty)$; accordingly, the jump function belongs to
$L^2(1,+\infty)$. On the other hand, more refined properties of continuity
and differentiability of the jump function are relevant in the mathematical
theory of the Cauchy integral equation, and particularly in the physical applications \cite{Bros2}.
In agreement with this approach,
we shall prove in Section \ref{sezione2} theorems which are variations on the Stein--Wainger
result. The main mathematical tool used in our approach is the Watson
resummation method which leads, in a very natural way, to the Laplace
transform of the jump function. It turns out that this Laplace transform coincides exactly
with the Carlsonian interpolation of the coefficients.
In this way all what is necessary for extending the methods
and the results to expansions in terms of Legendre and ultraspherical polynomials is obtained.
In this extension, in fact,
the interpolating function $\tilde{a}(\lambda)$ coincides with the
\emph{spherical Laplace transform} (in the sense of Faraut \cite{Bros1,Faraut1,Faraut2}),
that can be regarded as a composition of the classical Laplace
transform and the Abel--Radon transform. \\
This paper can be regarded as the completion and a large extension of a preliminary work by one of
the authors (G.A.V.) \cite{Viano}, where the Hausdorff moment problem has been
approached by the use of the Pollaczek polynomials.

\section{Hausdorff Moments, Hardy Spaces and Markov Processes}
\label{sezione1}
\subsection{Hausdorff Moments and Hardy Spaces}
\label{sezione1A}
Given a sequence of (real) numbers $\{f_n\}_0^\infty$, let $\Delta$ denote the difference operator:
\beq
\label{unouno}
\Delta f_n = f_{n+1} - f_n.
\eeq
Then, we have:
\beq
\label{unodue}
\Delta^k f_n = \underbrace{\Delta \times \Delta \times \cdots \times \Delta}_{k} f_n =
\sum_{m=0}^k (-1)^m \, \binom{\sds k}{\sds m} f_{n+k-m},
\eeq
(for every $k \geq 0$); $\Delta^0$ is the identity operator by definition.
Now, suppose that there exists a positive constant $M$ such that:
\beq
\label{unotre}
(n+1)\sum_{i=0}^n \binom{\sds n}{\sds i}^2 \left | \Delta^i f_{(n-i)} \right |^2 < M \qquad (n=0,1,2,\ldots).
\eeq
It can be proved \cite{Widder} that condition (\ref{unotre}) is necessary and sufficient
in order to represent the sequence $\left \{ f_n \right \}_0^\infty$ as follows:
\beq
\label{unoquattro}
f_n = \int_0^1 x^n u(x) \, dx \qquad (n=0,1,2,\ldots),
\eeq
where $u(x)$ belongs to $L^2(0,1)$.

We can prove the following Proposition.
\begin{proposition}
\label{pro:1}
If the sequence $\left \{ f_n \right \}_0^\infty$
satisfies condition $(\ref{unotre})$, then there exists a unique interpolation
of this sequence, denoted by
$\tilde{F}(\lambda)$ $(\lambda \in \C, [\tilde{F}(\lambda)]_{(\lambda=n)} = f_n)$
which belongs to the Hardy space $H^2(\C_{-1/2})$, and satisfies the following properties:
\begin{itemize}
\item[(i)] $\tilde{F}(\lambda)$ is holomorphic in the half--plane $\Real\lambda > -1/2$;
\item[(ii)] $\tilde{F}(\sigma+i\nu)$ belongs to $L^2(-\infty,+\infty)$ for
any fixed value of $\Real\lambda = \sigma \geq -1/2$;
\item[(iii)] $\tilde{F}(\lambda)$ tends uniformly to zero as $\lambda$
tends to infinity inside any fixed half--plane $\Real\lambda \geq \delta > -1/2$.
\end{itemize}
\end{proposition}

\begin{proof}
If the sequence $\{f_n\}_0^\infty$ satisfies condition (\ref{unotre}),
then representation (\ref{unoquattro}) holds true. If in this representation we
put $x=e^{-t}$, then we obtain:
\beq
\label{unocinque}
f_n = \int_0^{+\infty} e^{-nt} e^{-t} u(e^{-t}) \, dt \qquad (n=0,1,2,\ldots).
\eeq
Therefore the numbers $f_n$ can be regarded as the restriction to the integers
of the following Laplace transform:
\beq
\label{unosetteprimo}
\tilde{F}(\lambda) = \int_0^{+\infty} e^{-(\lambda+1/2)t} e^{-t/2} u(e^{-t}) \, dt.
\eeq
Indeed, one has $[\tilde{F}(\lambda)]_{(\lambda=n)}=f_n$. Moreover,
$\int_0^{+\infty}|\exp(-t/2)u(\exp(-t))|^2\,dt=\int_0^1|u(x)|^2\,dx<\infty$,
and therefore the function $\exp(-t/2)u(\exp(-t))$ belongs to $L^2(0,+\infty)$.
Then, in view of the Paley--Wiener theorem \cite{Hoffman} and of formula (\ref{unosetteprimo}),
we can conclude that $\tilde{F}(\lambda)$ belongs to the Hardy space $H^2(\C_{-1/2})$,
and properties (i), (ii) and (iii) follow.
Thus, we can make use of the Carlson theorem \cite{Boas}, which guarantees that
$\tilde{F}(\lambda)$ represents the unique interpolation of the sequence
$\left \{ f_n \right \}_0^\infty$.
\end{proof}

Now, we can prove the following Proposition.
\begin{proposition}
\label{pro:2}
If the sequence $\left \{ f_n \right \}_0^\infty$,
where $f_n = n^p a_n~(p \geq 1)$, satisfies condition $(\ref{unotre})$, then
there exists a unique Carlsonian interpolation of the numbers $\left \{ a_n \right \}_0^\infty$,
denoted by $\tilde{a}(\lambda)$ $(\lambda \in \C)$, which satisfies the following properties:
\begin{itemize}
\item[(i)] $\tilde{a}(\lambda)$ is holomorphic in $\Real\lambda >-1/2$;
\item[(ii)] $\lambda^p \tilde{a}(\lambda)$ belongs to $L^2(-\infty,+\infty)$ for any
fixed value of $\Real\lambda = \sigma \geq -1/2$;
\item[(iii)] $\lambda^p \tilde{a}(\lambda)$ tends uniformly to zero as $\lambda$
tends to infinity inside any fixed half--plane $\Real\lambda \geq \delta > -1/2$;
\item[(iv)] $\lambda^{(p-1)} \tilde{a}(\lambda)$ belongs to $L^1(-\infty,+\infty)$
for any fixed value of $\Real\lambda = \sigma \geq -1/2$.
\end{itemize}
\end{proposition}

\begin{proof}
Since the numbers $f_n = n^p a_n$ $(n=0,1,2,\ldots ; p \geq 1)$
satisfy condition (\ref{unotre}), then there exists a unique Carlsonian interpolation
of the sequence $\left \{ f_n \right \}_0^\infty$, denoted by $\tilde{F}(\lambda)$,
$(\lambda \in \C)$, which can be written as the product:
$\tilde{F}(\lambda) = \lambda^p \tilde{a}(\lambda)$, where $\tilde{a}(\lambda)$ is the unique
Carlsonian interpolation of the numbers $\left \{ a_n \right \}_0^\infty$.
In view of condition (\ref{unotre}) and Proposition \ref{pro:1} it follows that
$\tilde{F}(\lambda)$ belongs to $H^2(\C_{-1/2})$. Therefore properties
(i), (ii) and (iii) follow immediately. \\
By applying the Schwarz inequality, and recalling that $\tilde{F}(\lambda) \in L^2(-\infty,+\infty)$
for any fixed $\Real\lambda = \sigma \geq -1/2$, we have:
\beq
\label{unodieci}
\begin{split}
& \int_{-\infty}^{+\infty} \left |(\sigma+i\nu)^{(p-1)} \tilde{a}(\sigma+i\nu) \right |\,d\nu =
\int_{-\infty}^{+\infty} \left |\frac{\tilde{F}(\sigma+i\nu)}{(\sigma+i\nu)}\right |\,d\nu \\
& \qquad\leq \left ( \int_{-\infty}^{+\infty} \frac{1}{|(\sigma+i\nu)|^2}\,d\nu\right )^{1/2}
\left ( \int_{-\infty}^{+\infty} |\tilde{F}(\sigma+i\nu)|^2\, d\nu \right )^{1/2} < \infty,
\end{split}
\eeq
if $\sigma \geq -1/2$, $\sigma \neq 0$, $p \geq 1$. Finally, from inequality (\ref{unodieci}), and in view
of the regularity and integrability of the function $\lambda^{(p-1)}\tilde{a}(\lambda)$ in the
neighborhood of $\Real\lambda = 0$ we can state in all generality that $\lambda^{(p-1)}\tilde{a}(\lambda)$
belongs to $L^1(-\infty,+\infty)$ for any $p\geq 1$.
\end{proof}

\subsection{Hausdorff Moments and Markov Processes}
\label{sezione1B}
In this subsection we follow closely the paper of Watanabe \cite{Watanabe}. Let
$(\Omega, {\cF}, P)$ be an abstract probability field. If
$\{y_n(\tilde{\omega});\, n \geq 1 \}$ is a sequence of random variables
on $(\Omega, {\cF}, P)$ which are mutually independent, and each one  satisfies
\beq
\label{unoundici}
P\{y_n(\tilde{\omega})=1\}=p, \qquad P\{y_n(\tilde{\omega})=0\}=1-p,
\eeq
then it is called a Bernoulli sequence and denoted by $B(p)$. In the sequel we shall consider
$B(1/2)$. Let $E$ be the set of all points $(n,i)$ such that $n\geq i= 0,1,2,\cdots$.
Next, we consider the Markov process $x_n$ attached to $B(1/2)$. Let us note that:
\beq
\label{unododici}
P_{(n,i)}\{x_k=(m,j)\} =
\begin{cases}
\left (\sds\frac{1}{2}\right )^k \sds\binom{k}{j-i} & \text{for} \;\; m=n+k,\,j\geq i, \\
\; 0& \text{otherwise,}
\end{cases}
\eeq
where $k\geq 0$, $(n,i) \in E$, $(m,j) \in E$.
The kernel $K((n,i),(m,j))$ is given by \cite{Watanabe}:
\beq
\label{unotredici}
K((n,i),(m,j))= \frac{P_{(n,i)}\{\sigma(\{m,j\})<+\infty\}}{P_{(0,0)}\{\sigma(\{m,j\})<+\infty\}}=
2^n \frac{(m-n)!\, j!\, (m-j)!}{m!\, (m-n-j+i)!\, (j-i)!}.
\eeq
Now, consider an infinite sequence $(m_k,j_k)$ having no limit point in $E$, and
such that
\beq
\label{unoquattordici}
\lim_{k\rightarrow\infty}\frac{j_k}{m_k}=1-b,
\eeq
for a suitable $0\leq b\leq 1$. Then, using the Stirling formula, and taking into
account equality (\ref{unoquattordici}), from (\ref{unotredici}) we obtain:
\beq
\label{unoquindici}
K((n,i),b)=2^n b^{(n-i)}(1-b)^i.
\eeq
Thus, one may consider that the Martin boundary ${\cM}$ \cite{Hunt} induced
by the process $x_n$ coincides with the interval $[0,1]$ as a set, and,
accordingly, the generalized Poisson kernel $K((n,i),b)$ is given by
$2^n b^{(n-i)}(1-b)^i$. Finally, we note that for a function $u$ over $E$
the expectation $E_{(0,0)}$ is given by:
\beq
\label{unosedici}
E_{(0,0)}(|u(x_n)|)=2^{-n}\sum_{i=0}^n |u(n,i)|\binom{\sds n}{\sds i}.
\eeq
Now, the following propositions due to Watanabe \cite{Watanabe} can be stated.
\begin{proposition}
\label{pro:3}
Let $x_n$ be the Markov process attached to
the Bernoulli sequence $B(1/2)$.
\begin{itemize}
\item[(i)] The Martin boundary induced by $x_n$ is equivalent to the
interval $[0,1]$ with the ordinary topology;
\item[(ii)] The generalized Poisson kernel $K((n,i),b)$ is:
\beq
\label{dueduetre}
K((n,i),b) = 2^n b^{(n-i)}(1-b)^i.
\eeq
\item[(iii)] A function $u$ (belonging to the set of all the finite real
valued functions over $E^\star=E \cup \infty$, vanishing at $\infty$) can
be represented by means of a bounded signed measure on
$\left ( [0,1], {\cB}_{[0,1]}\right )$, (where ${\cB}_{[0,1]}$ is the
Borel field consisting of all the ordinary Borel subsets in $[0,1]$), as follows:
\beq
\label{unodiciassette}
u(n,i)=2^n\int_0^1 b^{(n-i)} (1-b)^i \, d\mu(b),
\eeq
(for every $(n,i)\in E$), if and only if $u$ is $x_n$--harmonic, and
$E_{(0,0)}(|u(x_n)|)$ is bounded in $n$.
\end{itemize}
\end{proposition}

\begin{proof}
See \cite{Watanabe}.
\end{proof}

\begin{proposition}
\label{pro:4}
Let $x_n$ be the Markov process attached to
the Bernoulli sequence $B(1/2)$. Given a sequence of real numbers
$\{f_n;\,n\geq 0\}$ such that:
\beq
\label{unodiciotto}
\sum_{i=0}^n \left |\Delta^i f_{(n-i)}\right |\binom{\sds n}{\sds i} < L \qquad (n=0,1,2,\ldots; L=\text{\rm constant}),
\eeq
then the function $u(n,i)$ defined by:
\beq
\label{unodiciannove}
u(n,i)=2^n(-1)^i \Delta^i f_{(n-i)}
\eeq
is a $x_n$--harmonic function, and can be represented by formula
(\ref{unodiciassette}).
\end{proposition}

\begin{proof}
See \cite{Watanabe}.
\end{proof}

Notice that from representation (\ref{unodiciassette}) it follows:
\beq
\label{unoventi}
2^{-n} u(n,0)=\int_0^1 b^n \,d\mu(b),
\eeq
which can be compared with representation (\ref{unoquattro}).
Moreover, if the sequence $\{f_n\}_0^\infty$ satisfies
inequality (\ref{unotre}), then it satisfies also inequality (\ref{unodiciotto}).
This can be proved easily by the use of the Cauchy inequality:
\beq
\label{unoventuno}
\sum_{i=0}^n |x_i \, y_i | \leq \left (\sum_{i=0}^n |x_i|^2\right )^{1/2}
\left (\sum_{i=0}^n |y_i|^2\right )^{1/2}.
\eeq
In fact, if in inequality (\ref{unoventuno}) we put: $y_i=1$, $\forall i \in (0,1,2,\ldots,n)$,
$x_i=|\Delta^i f_{(n-i)}|\binom{n}{i}$, we obtain:
\beq
\label{unoventidue}
\sum_{i=0}^n \binom{\sds n}{\sds i} |\Delta^i f_{(n-i)}| \leq
\left\{\sum_{i=1}^n\binom{\sds n}{\sds i}^2 |\Delta^i f_{(n-i)}|^2\right\}^{1/2}
(n+1)^{1/2}.
\eeq
Therefore, from inequalities (\ref{unotre}) and (\ref{unoventidue}) we obtain:
\beq
\label{unoventitre}
\sum_{i=0}^n \binom{\sds n}{\sds i} |\Delta^i f_{(n-i)}| \leq
\frac{\sqrt{M}}{(n+1)^{1/2}} (n+1)^{1/2} = \sqrt{M},
\eeq
that coincides with inequality (\ref{unodiciotto}), if we
put $L=\sqrt{M}$.

\section{A Double Analytic Structure for a Class of Trigonometric and Power Series}
\label{sezione2}
In the complex plane $\C$ of the variable $\theta=u+iv$ ($u,v \in \R$) we
consider the following domains:
${\cI}_+^{(\pm\xi_0)} = \{\theta \in \C \mid \Imag\theta > \pm \xi_0,\,\xi_0 \geq 0\}$,
and
${\cI}_-^{(\pm\xi_0)} = \{\theta \in \C \mid \Imag\theta < \pm \xi_0,\,\xi_0 \geq 0\}$.
We introduce, correspondingly, the following \emph{cut domains}:
${\cI}_+^{(0)}\setminus \Xi_+^{(\xi_0)}$, where
$\Xi_+^{(\xi_0)}=\{\theta\in\C \mid\theta=2k\pi + iv,\, v>\xi_0,\, \xi_0\geq 0,\, k\in \Z\}$,
and
${\cI}_-^{(0)}\setminus \Xi_-^{(-\xi_0)}$, where
$\Xi_-^{(-\xi_0)}=\{\theta\in\C \mid\theta=2k\pi + iv,\, v<-\xi_0,\, \xi_0\geq 0,\, k\in \Z\}$
(for a detailed description of these cut domains see \cite{Bros1}).
We will use the notation $\dot{A}=A\setminus 2\pi\Z$ for every subset $A$ of $\C$
which is invariant under the translation group $2\pi\Z$. \\
We can then prove the following theorem.

\begin{figure}[ht]
\begin{center}
\includegraphics[width=\textwidth]{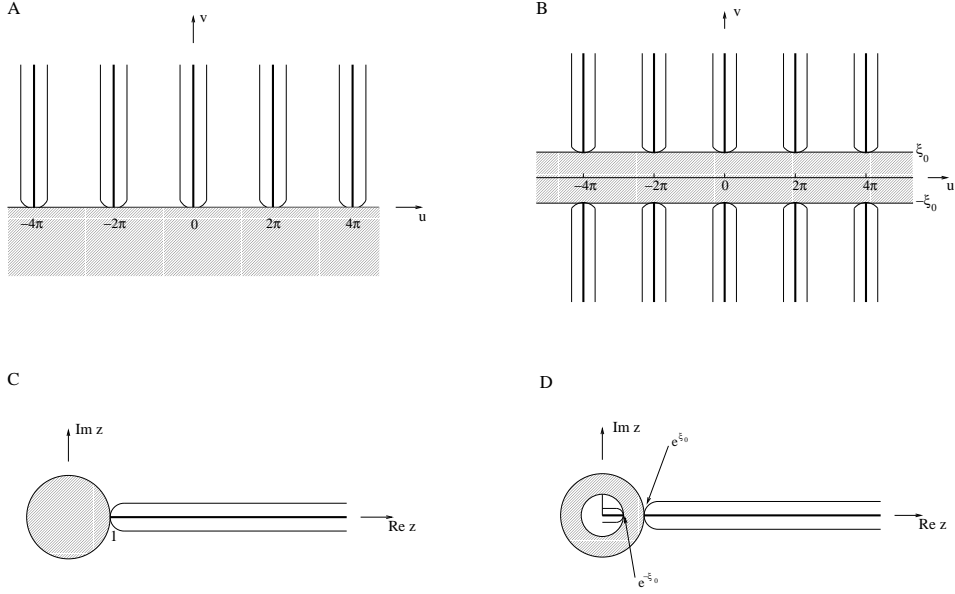}
\caption{\label{figura_1}
\emph{Cut--domains} in the various geometries. The dashed regions represent
the analyticity domains of the functions $f$ to which the series converges, whereas
the thick lines indicate the cuts that delimit the holomorphic extension of $f$.
A: See Theorem \protect\ref{the:1}.
B: See Proposition \protect\ref{pro:5}.
C: See Theorem \protect\ref{the:1prime}.
D: See Proposition \protect\ref{pro:6}.}
\end{center}
\end{figure}

\begin{theorem}
\label{the:1}
Let us consider the following series:
\beq
\label{dueuno}
\frac{1}{2\pi}\sum_{n=0}^\infty a_n e^{-in\theta} \qquad (\theta=u+iv,\, u,v\in\R),
\eeq
and suppose that the set of numbers $\{f_n\}_0^\infty$, $f_n=n^pa_n$ $(p\geq 1,\, n=0,1,2,\ldots)$
satisfies condition $(\ref{unotre})$, then:
\begin{itemize}
\item[(1)] series $(\ref{dueuno})$ converges uniformly to a function $f(\theta)$
analytic in ${\cI}_-^{(0)}$;
\item[(2)] the function $f(\theta)$ admits a holomorphic extension to the
\emph{cut domain} ${\cI}_+^{(0)}\setminus \dot{\Xi}_+^{(0)}$ (see Fig. $\ref{figura_1}$A);
\item[(3)] the jump function $F(v)$ (which equals the discontinuity of $if(\theta)$
across the cuts $\dot{\Xi}_+^{(0)}$) is a function of class $C^{(p-1)}$ ($p\geq 1$),
and satisfies the following bound:
\beq
\label{duedue}
|F(v)| \leq \|\tilde{a}_\sigma\|_1 \, e^{\sigma v} \qquad (\sigma \geq -\sfrac{1}{2},\, v\in\R^+),
\eeq
where $\tilde{a}(\sigma+i\nu)$ $(\nu\in\R)$ is the Carlsonian interpolation of the
coefficients $a_n$, and
\beq
\label{duetre}
\|\tilde{a}_\sigma\|_1 = \frac{1}{2\pi}\int_{-\infty}^{+\infty} \left |\tilde{a}(\sigma+i\nu) \right |\,
d\nu \qquad (\sigma \geq -\sfrac{1}{2});
\eeq
\item[(4)] $\tilde{a}(\sigma+i\nu)$ is the Laplace transform of the jump function $F(v)$: i.e.,
\beq
\label{duequattro}
\tilde{a}(\sigma+i\nu) = \int_0^{+\infty} F(v) e^{-(\sigma+i\nu)v}\, dv \qquad (\sigma > -\sfrac{1}{2});
\eeq
\item[(5)] the Plancherel equality holds true:
\beq
\label{duecinque}
\int_{-\infty}^{+\infty} |\tilde{a}(\sigma+i\nu)|^2\,d\nu = 2\pi\int_0^{+\infty} |F(v) e^{-\sigma v}|^2\,dv
 \qquad (\sigma \geq -\sfrac{1}{2}).
\eeq
\end{itemize}
\end{theorem}

\begin{proof}
Since the set $\{f_n\}_0^\infty$ satisfies condition (\ref{unotre}),
given an arbitrary number $C$, there exists a real number $m$ such that for $n>m$,
$|a_n| \leq C$. Therefore, we can write:
\beq
\label{duesei}
\left |\sum_{n>m}^\infty \left (\frac{1}{2\pi} a_n e^{-in\theta}\right )\right | \leq
\frac{C}{2\pi} \sum_{n>m}^\infty e^{nv} \qquad (\theta=u+iv).
\eeq
The series at the r.h.s. of formula (\ref{duesei}) is uniformly convergent for $v\leq v_0<0$.
Recalling the Weierstrass theorem on the uniformly convergent series of
analytic functions, we can also conclude that the series
$(1/2\pi)\sum_{n>m}^\infty a_n\exp(-in\theta)$ $(\theta=u+iv)$ converges uniformly
to a function analytic in ${\cI}_-^{(0)}$. On the other hand, series
(\ref{dueuno}) can be rewritten as the following sum:
\beq
\label{duesette}
\frac{1}{2\pi} \sum_{n=0}^\infty a_n e^{-in\theta}=
\frac{1}{2\pi}\left\{\sum_{n>m}^\infty a_n e^{-in\theta}+T_m(\theta)\right\},
\eeq
where $T_m(\theta)=\sum_{n=0}^{[m]}(a_n e^{-in\theta})$ is a trigonometric
polynomial analytic in ${\cI}_-^{(0)}$. Therefore the first statement is proved.

In order to prove the other statements, let us introduce the following integral:
\beq
\label{dueotto}
f_\epsilon(u)=\frac{i}{4\pi}\int_{{\cC}}\frac{\tilde{a}(\lambda) e^{-i\lambda(u-\epsilon\pi)}}
{\sin\pi\lambda}\, d\lambda \qquad (\epsilon=\pm),
\eeq
where $\tilde{a}(\lambda)$, $(\lambda=\sigma+i\nu)$ is the unique Carlsonian interpolation
of the sequence $\{a_n\}_0^\infty$, which exists in view of the fact that the set
$\{f_n\}_0^\infty$ satisfies condition (\ref{unotre}), and
the contour ${\cC}$ is contained in the half--plane
$\C_{-1/2}$
and encircles the positive real semi--axis of the $\lambda$--plane
(or a part of it) as is illustrated in Fig. \ref{figura_2}A.

\begin{figure}[ht]
\begin{center}
\includegraphics[width=7cm]{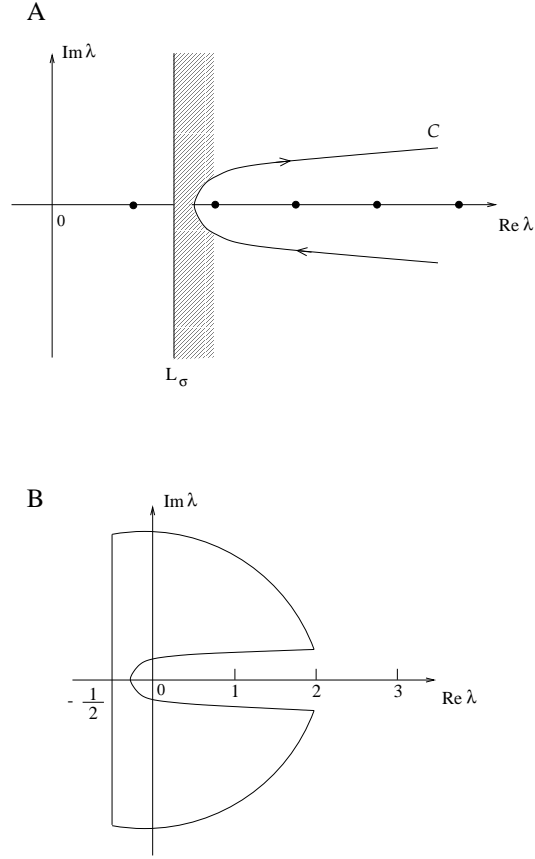}
\caption{\label{figura_2}
A: ${\cC}$ is the integration path of the integral (\protect\ref{dueotto}). The straight
line path $L_\sigma$ is here represented for $\sigma = 3/2$ (see Theorem \protect\ref{the:1}).
B: Contour integration used for evaluating integral (\protect\ref{trenove}) (see Theorem \protect\ref{the:2}).}
\end{center}
\end{figure}

Now, let us consider the following inequalities:
\begin{eqnarray}
\label{duenove}
&&\left |e^{-i(\sigma+i\nu)(u-\pi)}\right | \leq 2\cosh\pi\nu~ \qquad (u\in [0,2\pi]), \\
\label{duedieci}
&&|\sin\pi (\sigma+i\nu)| \geq \sinh\pi\nu, \\
\label{dueundici}
&&\left | \frac{e^{-i(\sigma+i\nu)(u-\pi)}}{\sin\pi(\sigma+i\nu)}\right | \leq
\left | \frac{2\cosh\pi\nu}{\sinh\pi\nu}\right | \qquad (u\in [0,2\pi]).
\end{eqnarray}
Let us recall that $\lambda^p\tilde{a}(\lambda)$ ($p \geq 1$) tends uniformly to zero as
$\lambda\rightarrow\infty$ inside any fixed half--plane $\Real\lambda=\sigma\geq\delta >-1/2$,
and $\lambda^{(p-1)}\tilde{a}(\lambda)$ belongs to $L^1(-\infty,+\infty)$ for
$\Real\lambda=\sigma \geq -1/2$ (see Proposition \ref{pro:2}).

\noindent
In view of these properties of $\tilde{a}(\lambda)$, and by the use of bound (\ref{dueundici}),
we can guarantee that the integral $f_+(u)$ $(u\in [0,2\pi])$ converges, and the contour
${\cC}$ can be deformed and replaced by the line
$L_\sigma=\{\lambda=\sigma+i\nu,\, \nu\in\R , \sigma\geq -1/2\}$, provided that the
real variable $u$ is kept in $[0,2\pi]$ (see Fig. \ref{figura_2}A). Finally, by applying
the Watson resummation method \cite{Watson} we obtain for $u\in [0,2\pi]$:
\beq
\label{duedodici}
f_+(u) = -\frac{1}{4\pi}\int_{-\infty}^{+\infty} \frac{ \tilde{a}(\sigma+i\nu) e^{-i(\sigma+i\nu)(u-\pi)} }
{\sin\pi (\sigma+i\nu)} \,d\nu =
\frac{1}{2\pi}\sum_{n=l}^\infty a_n e^{-inu},
\eeq
($l \geq 0$, integer; $-1/2\leq\sigma < 0$ if $l=0$, $l-1<\sigma <l$ if $l>0$).

Proceeding in an analogous fashion for the integral $f_-(u)$ (formula (\ref{dueotto}) with
$\epsilon = -$), and distorting the contour integration in a similar way, we finally
obtain for $u\in [-2\pi,0]$
\beq
\label{duetredici}
f_-(u) = -\frac{1}{4\pi}\int_{-\infty}^{+\infty} \frac{ \tilde{a}(\sigma+i\nu) e^{-i(\sigma+i\nu)(u+\pi)} }
{\sin\pi (\sigma+i\nu)} \,d\nu =
\frac{1}{2\pi}\sum_{n=l}^\infty a_n e^{-inu},
\eeq
($l \geq 0$, integer; $-1/2\leq\sigma < 0$ if $l=0$, $l-1<\sigma <l$ if $l>0$).

Now, in the integral (\ref{duedodici}) we substitute for $u$ the complex variable
$\theta=u+iv$, and we see that the obtained integral provides an analytic continuation of $f_+$
in the strip $\{\theta=u+iv,\,0<u<2\pi\}$, continuous in the closure of the latter. Indeed we have:
\beq
\label{duequattordici}
e^{-\sigma v} f_+(u+iv)=\frac{1}{2\pi}\int_{-\infty}^{+\infty}e^{i\nu v}H_\sigma^u(\nu)\,d\nu
 \qquad (0\leq u\leq 2\pi),
\eeq
with
\beq
\label{duequindici}
H_\sigma^u(\nu)=-\frac{\tilde{a}(\sigma+i\nu)e^{-i(\sigma+i\nu)(u-\pi)}} {2\sin\pi(\sigma+i\nu)},
\eeq
then, in view of bound (\ref{dueundici}), and since
$\lambda^{(p-1)}\tilde{a}(\lambda) \in L^1(-\infty,+\infty)$ for any fixed value of
$\Real\lambda=\sigma\geq -1/2$ (see Proposition \ref{pro:2}), the statement above is proved.
Similarly, the analytic continuation of the function $f_-$ is defined in the strip
$\{\theta=u+iv~,\,-2\pi < u < 0\}$. The discontinuity
$f_+(iv)-f_-(iv)$ can be computed by replacing $u$ by $iv$ in integrals (\ref{duedodici}) and
(\ref{duetredici}), and subtracting Eq. (\ref{duetredici}) from Eq. (\ref{duedodici}).
We then obtain
\beq
\label{duesedici}
i[f_+(iv)-f_-(iv)]=
\frac{1}{2\pi}\int_{-\infty}^{+\infty}\tilde{a}(\sigma +i\nu)e^{(\sigma+i\nu)v}\,d\nu \qquad (v\in\R^+,\,\sigma\geq -\sfrac{1}{2}).
\eeq
Thus, we have proved that the function $f(\theta)$ $(\theta=u+iv)$ admits a holomorphic extension
to the \emph{cut domain} ${\cI}_+^{(0)}\setminus\dot{\Xi}_+^{(0)}$.

From formula (\ref{duesedici}) we derive the following bound for the jump function
$F(v) :\equiv i[f_+(iv)-f_-(iv)]$:
\beq
\label{duediciassette}
|F(v)| \leq \|\tilde{a}_\sigma\|_1 e^{\sigma v} \qquad (\sigma\geq -\sfrac{1}{2},\, v\in\R^+),
\eeq
where
\beq
\label{duediciotto}
\|\tilde{a}_\sigma\|_1 = \frac{1}{2\pi}\int_{-\infty}^{+\infty} |\tilde{a}(\sigma+i\nu)|\,d\nu < \infty
\qquad (\sigma\geq -\sfrac{1}{2}).
\eeq
Using again formula (\ref{duesedici}), and recalling the Riemann--Lebesgue theorem, we can prove
that $F(v)$ is a function of class $C^{(p-1)}$ $(p\geq 1)$ in view of the
fact that $\lambda^{(p-1)}\tilde{a}(\lambda)$ belongs to $L^1(-\infty,+\infty)$ for
any $\Real\lambda=\sigma\geq -1/2$ (see Proposition \ref{pro:2}).

Inverting formula (\ref{duesedici}) we obtain:
\beq
\label{duediciannove}
\tilde{a}(\sigma+i\nu)=\int_0^{+\infty}F(v)e^{-(\sigma+i\nu)v}\,dv \qquad (\sigma > -\sfrac{1}{2}),
\eeq
which is, indeed, the Laplace transform of the jump function $F(v)$, and it is holomorphic
for $\Real\lambda=\sigma>-1/2$.

\noindent
Finally, recalling that $\tilde{a}(\sigma+i\nu)$ belongs to $L^2(-\infty,+\infty)$ at any
fixed value $\Real\lambda=\sigma\geq -1/2$, we obtain the Plancherel equality (\ref{duecinque})
and, in particular:
\beq
\label{duediciannoveprimo}
\int_{-\infty}^{+\infty}\left |\tilde{a}\left(-\frac{1}{2}+i\nu\right) \right |^2\,d\nu=
2\pi\int_0^{+\infty}|F(v)e^{v/2}|^2\,dv.
\eeq
\end{proof}

\begin{remarks}
(i) Theorem \ref{the:1} proves a double analytic structure, connected with
series (\ref{dueuno}), in the following sense: To the functions $\tilde{a}(\lambda)$,
that interpolate the coefficients $a_n$, and are holomorphic in the half--plane
$\Real\lambda >-1/2$, there corresponds the class of functions $f(\theta)$ holomorphic
in the domain ${\cI}_-^{(0)} \cup ({\cI}_+^{(0)}\setminus \dot{\Xi}_+^{(0)})$, and moreover
$\tilde{a}(\lambda)$ is the Laplace transform of the jump function $F(v)$ across the cut. \\
(ii) Let us note that the Plancherel equality (formulae (\ref{duecinque})
and (\ref{duediciannoveprimo})), as well as the analyticity property of the function $f(\theta)$,
remain true under milder conditions on the coefficients $a_n$. In fact, it is sufficient that
the $a_n$'s form a sequence of numbers that satisfies condition (\ref{unotre}).
This is, indeed, the result contained in the theorem of Stein--Wainger \cite{Stein} referred to
series (\ref{dueuno}). In conclusion, the more restrictive conditions, assumed in Theorem
\ref{the:1}, are reflected by the smoothness property of the jump function, which is, however, a quite relevant
property playing an important role in the applications to physical problems.
\end{remarks}

By substituting the complex plane of the variable $z=\exp(-i\theta)$ to
the $2\pi$--periodic $\theta$--plane,
we can now give an equivalent presentation of the results of Theorem \ref{the:1}, in terms of
properties of Taylor series and of Mellin transformation. To the cut at $\theta=iv$, it corresponds,
in the $z$--plane geometry, the cut located on the real axis ($x \equiv \Real z$) from
$x=1$ up to $+\infty$. To the jump function $F(v)$ it corresponds the function
$F(\ln x)$, which shall still be denoted hereafter simply by $F(x)$ with a small abuse of
notation which avoids, however, an useless proliferation of symbols. Adopting the same convention,
we shall always denote the jump function, in the various geometries, with the same
symbol: i.e., $F$.

\begin{theorem}
\label{the:1prime}
If in the Taylor series:
\beq
\label{dueventi}
\frac{1}{2\pi}\sum_{n=0}^\infty a_n z^n \qquad (z=x+iy;\,x,y\in\R),
\eeq
the coefficients $a_n$ satisfy the assumptions required by Theorem $\ref{the:1}$, then:
\begin{itemize}
\item[(1')] the series converges uniformly to a function $f(z)$ analytic in the unit disk
$D_0=\{z \mid |z|<1\}$;
\item[(2')] $f(z)$ admits a holomorphic extension to the \emph{cut plane} $z\in\C\setminus (1,+\infty)$
(see Fig. $\ref{figura_1}$C);
\item[(3')] the jump function $F(x)=-i(f_+(x)-f_-(x))$,
($f_\pm(x)=\lim_{\stackrel{\scriptstyle\epsilon\rightarrow 0}{\epsilon >0}}
f(x\pm i\epsilon$)) is a function of class $C^{(p-1)}$ $(p\geq 1)$, and satisfies
the following bound:
\beq
\label{dueventuno}
|F(x)|\leq\|\tilde{a}_\sigma\|_1 x^\sigma \qquad (\sigma\geq -\sfrac{1}{2},\,x\in(1,+\infty)),
\eeq
where $\tilde{a}(\sigma+i\nu)$ $(\nu\in\R)$ is the Carlsonian interpolation of the
coefficients $a_n$, and $\|\tilde{a}_\sigma\|_1$ is given by formula $(\ref{duetre})$;
\item[(4')] $\tilde{a}(\sigma+i\nu)$ is the Mellin transform of the jump function: i.e.,
\beq
\label{dueventidue}
\tilde{a}(\sigma+i\nu)=\int_1^{+\infty}F(x)x^{-(\sigma+i\nu)-1}\,dx \qquad (\sigma>-\sfrac{1}{2});
\eeq
\item[(5')] the Plancherel formula associated with the Mellin transform gives
\beq
\label{dueventitre}
\int_{-\infty}^{+\infty}|\tilde{a}(\sigma+i\nu)|^2\,d\nu=
2\pi\int_1^{+\infty}|F(x)|^2 x^{-2\sigma-1}\,dx \qquad (\sigma\geq -\sfrac{1}{2}).
\eeq
\end{itemize}
\end{theorem}

\begin{remark}
Equality (\ref{dueventitre}) that, in the particular case
$\sigma=-1/2$, reads
\beq
\label{dueventiquattro}
\int_{-\infty}^{+\infty}\left |\tilde{a}\left(-\frac{1}{2}+i\nu\right)\right |^2\,d\nu=2\pi\int_1^{+\infty}|F(x)|^2\,dx,
\eeq
coincides with the result contained in the Stein--Wainger theorem \cite{Stein}. In order to obtain
this latter result it is sufficient to require that the set $\{a_n\}_0^\infty$ satisfies
condition (\ref{unotre}).
\end{remark}

We now present, without giving the proof, two variants of Theorems \ref{the:1} and \ref{the:1prime}, which are
relevant in the physical applications, and specifically in the theory of the thermal Green
functions \cite{Bros2}. The proofs can be easily obtained, with small variations, from those
of Theorems \ref{the:1} and \ref{the:1prime}.
\begin{proposition}
\label{pro:5}
If in the following series
\beq
\label{dueventisei}
\frac{1}{2\pi}\sum_{n=-\infty}^{+\infty} a_n \cos nu \qquad (u\in \R),
\eeq
the coefficients $a_n$, in addition to the assumptions required by Theorem $\ref{the:1}$,
satisfy also the following bound:
\beq
\label{dueventicinque}
|a_n| \leq C \, e^{-(n-m)\xi_0} \qquad (n>m,\,m\in\R^+,\xi_0>0),
\eeq
then
\begin{itemize}
\item[(1)] series (\ref{dueventisei}) converges uniformly to a function $f(\theta)$
analytic in the strip $\{\theta\in\C \mid |\Imag\theta | < \xi_0\}$;
\item[(2)] the function $f(\theta)$ admits a holomorphic extension to the
\emph{cut domain}
${\cI}_+^{(0)}\setminus\dot{\Xi}_+^{(\xi_0)} \cup {\cI}_-^{(0)}\setminus\dot{\Xi}_-^{(-\xi_0)}$
(see Fig. $\ref{figura_1}$B);
\item[(3)] the jump function across the cuts $\dot{\Xi}_\pm^{(\pm\xi_0)}$ satisfies
all the properties proved in Theorem $\ref{the:1}$.
\end{itemize}
\end{proposition}

\begin{proposition}
\label{pro:6}
If in the following Laurent series
\beq
\label{duetrentadue}
\frac{1}{2\pi}\sum_{n=-\infty}^{+\infty}a_n z^n \qquad (z=x+iy,\,x,\,y\in\R ,\, a_n=a_{-n}),
\eeq
the coefficients $a_n$ satisfy the assumption of Theorem $\ref{the:1prime}$ and condition $(\ref{dueventicinque})$,
then:
\begin{itemize}
\item[(1)] series $(\ref{duetrentadue})$ converges uniformly to a function $f(z)$
analytic in the annulus ${\cA}=\{z\in\C \mid \exp(-\xi_0) < |z| < \exp(\xi_0)\}$;
\item[(2)] $f(z)$ admits a holomorphic extension to the \emph{cut domain}
$z\in\C\setminus [(0,\exp(-\xi_0) \cup (\exp(\xi_0),+\infty)]$ (see Fig. $\ref{figura_1}$D);
\item[(3)] the jump function $F(x)$ satisfies the following bounds:
\beq
\label{duetrentatre}
|F(x)|\leq\|\tilde{a}_\sigma\|_1 x^\sigma \qquad (\sigma\geq -\sfrac{1}{2},\, x\in (e^{\xi_0},+\infty)),
\eeq
and, accordingly,
\beq
\label{duetrentaquattro}
\left |F\left (\frac{1}{x}\right )\right |\leq\|\tilde{a}_\sigma\|_1 x^{-\sigma} \qquad
(\sigma\geq -\sfrac{1}{2},\, x\in (0, e^{-\xi_0})),
\eeq
where $\tilde{a}(\sigma+i\nu)$ is the Carlsonian interpolation of the coefficients $a_n$,
and $\|\tilde{a}_\sigma\|_1$ is given by formula $(\ref{duetre})$;
\item[(4)] $\tilde{a}(\sigma+i\nu)$ is the Mellin transform of the jump function, and it is given by:
\beq
\label{duetrentacinque}
\tilde{a}(\sigma+i\nu)=\int_1^{+\infty}F(x) x^{-(\sigma+i\nu)-1}\, dx \qquad (\sigma>-\sfrac{1}{2}),
\eeq
or, equivalently, by:
\beq
\label{duetrentasei}
\tilde{a}(\sigma+i\nu)=\int_0^1 F\left (\frac{1}{x}\right ) x^{(\sigma+i\nu)-1}\, dx \qquad (\sigma>-\sfrac{1}{2}).
\eeq
\item[(5)] The Plancherel formula associated with the Mellin transform gives:
\beq
\label{duetrentasette}
\int_{-\infty}^{+\infty} |\tilde{a}(\sigma+i\nu)|^2\, d\nu =
2\pi\int_1^{+\infty} |F(x)|^2 x^{-2\sigma -1}\, dx \qquad (\sigma\geq -\sfrac{1}{2}),
\eeq
and, accordingly,
\beq
\label{duetrentotto}
\int_{-\infty}^{+\infty} |\tilde{a}(\sigma+i\nu)|^2\, d\nu =
2\pi\int_0^1 \left |F\left (\frac{1}{x}\right )\right |^2 x^{(2\sigma -1)}\, dx \qquad (\sigma\geq - \sfrac{1}{2}).
\eeq
\end{itemize}
\end{proposition}

\begin{remark}
Let us suppose that instead of the condition assumed in Theorem \ref{the:1}, the coefficients $a_n$
of series (\ref{dueuno}) satisfy the following condition: The sequence
$\{f_n\}_{n_0}^\infty$ $(n_0>0,\,f_n=n^p a_n,\, p\geq 1)$ satisfies condition (\ref{unotre}).
In this case all the results proved above remain true, except that now the constant $\sigma$,
that controls the bounds on the jump function, is larger than or equal to
$(n_0-1/2)$ $(n_0>0)$. Accordingly, the Laplace transform of the jump function (see formula
(\ref{duequattro})) holds true only for $\sigma >(n_0-1/2)$, and the Plancherel equality
(\ref{duecinque}) for $\sigma\geq (n_0-1/2)$ $(n_0>0)$. Analogous modifications should be
considered in relation to the other expansions treated in Theorem \ref{the:1prime} and in Propositions \ref{pro:5} and \ref{pro:6}.
\end{remark}

\section{Solution of a Class of Cauchy Integral Equations:
Representation of the Jump Function in Terms of the Coefficients $\bf a_n$}
\label{sezione3}
Let us now focus our attention on the series (\ref{dueuno}); the results which we obtain can be
easily extended to the other series that have been considered in Section \ref{sezione2}.
Here we assume on the
coefficients $a_n$ the weakest possible condition: i.e., the sequence $\{a_n\}_0^\infty$ is supposed
to satisfy condition (\ref{unotre}). Under this condition we can still guarantee that there exists
a unique Carlsonian interpolation $\tilde{a}(\lambda)$ of the coefficients $a_n$, and also that
$\tilde{a}(\sigma+i\nu)$ belongs to $L^2(-\infty,+\infty)$ for any fixed
$\Real\lambda=\sigma\geq -1/2$; but $\tilde{a}(-1/2+i\nu)$ does not belong, in general, to
$L^1(-\infty,+\infty)$. Therefore the inversion of the Fourier transform at $\sigma=-1/2$
holds only as a limit in the mean order two, which reads:
\beq
\label{treuno}
F(v)e^{v/2}=~ \lm\displaylimits_{\nu_0\rightarrow +\infty} \left (\frac{1}{2\pi}\int_{-\nu_0}^{\nu_0}
\tilde{a}\left(-\frac{1}{2}+i\nu\right)e^{i\nu v}\,d\nu\right ) \qquad (v\in\R^+).
\eeq
We can prove the following theorem.

\begin{theorem}
\label{the:2}
If in series $(\ref{dueuno})$ the coefficients $a_n$ satisfy
condition $(\ref{unotre})$, then the function $\exp(v/2)F(v)$ can be represented by the following expansion,
that converges in the sense of the $L^2$--norm:
\beq
\label{tredue}
e^{v/2}F(v) = \sum_{m=0}^\infty c_m\Phi_m(v) \qquad (v\in\R^+),
\eeq
where
\beq
\label{tretre}
c_m=\sqrt{2}\,\sum_{n=0}^\infty \frac{(-1)^n}{n!}\, a_n \,P_m\left [-i\left ( n+\frac{1}{2}\right )\right ],
\eeq
and
\beq
\label{trequattro}
\Phi_m(v)=i^m\,\sqrt{2} \,L_m(2e^{-v}) \,e^{{-e}^{-v}} e^{-v/2},
\eeq
$P_m$ and $L_m$ being, respectively, the Pollaczek and the Laguerre polynomials.
\end{theorem}

\begin{proof}
The Pollaczek polynomials $P_m^{(\alpha)}(\nu)$ $(\nu\in\R)$ are a set of polynomials
orthogonal in $(-\infty,+\infty)$ with the weight function (see \cite{Bateman,Szego}):
\beq
\label{trecinque}
w(\nu)=\frac{1}{\pi}\, 2^{(2\alpha -1)} |\Gamma(\alpha+i\nu)|^2 \qquad (\alpha > 0),
\eeq
(where $\Gamma(x)$ denotes the Euler gamma function).\\
We put $\alpha=1/2$ (in the following we shall omit the index in the notation of the
Pollaczek polynomials). Then the property of orthogonality reads as follows:
\beq
\label{tresei}
\int_{-\infty}^{+\infty} w(\nu)P_k(\nu)P_l(\nu)\, d\nu = \delta_{k,l},
\eeq
where now $w(\nu)=(1/\pi)|\Gamma(1/2+i\nu)|^2$.

\noindent
Next, we introduce the following functions (that may be called Pollaczek functions):
\beq
\label{tresette}
\psi_m(\nu)=\frac{1}{\sqrt{\pi}}\Gamma\left (\frac{1}{2}+i\nu\right )P_m(\nu) \qquad (P_0(\nu)=1),
\eeq
which form a complete basis in $L^2(-\infty,+\infty)$ (see \cite{Itzykson}).

\noindent
In view of the fact that the coefficients $\{a_n\}_0^\infty$ satisfy condition (\ref{unotre}),
then there exists a unique Carlsonian interpolation of the sequence $\{a_n\}_0^\infty$, denoted
by $\tilde{a}(\lambda)$, such that $\tilde{a}(-1/2+i\nu)$ $(\nu\in\R)$ belongs to $L^2(-\infty,+\infty)$
(see Proposition \ref{pro:1}). Therefore the function $\tilde{a}(-1/2+i\nu)$ can be expanded in terms of
the Pollaczek functions as follows:
\beq
\label{treotto}
\tilde{a}\left (-\frac{1}{2}+i\nu\right )=\sum_{m=0}^\infty d_m \psi_m(\nu) \qquad (\nu\in\R),
\eeq
and the convergence of this expansion is in the sense of the $L^2$--norm.

\noindent
The coefficients $d_m$ are given by:
\beq
\label{trenove}
d_m=\frac{1}{\sqrt{\pi}}\int_{-\infty}^{\infty} \tilde{a}\left (-\frac{1}{2}+i\nu\right )\,
\Gamma\left (\frac{1}{2}-i\nu\right )
\, P_m(\nu)\, d\nu.
\eeq
Taking into account the asymptotic behavior of the gamma function, we may evaluate integral
(\ref{trenove}) by the contour integration method along the path shown in Fig. \ref{figura_2}B.
Note that the poles of the gamma function $\Gamma(1/2-i\nu)$ are located at $\nu=-i(n+1/2)$,
and $[\tilde{a}(-1/2+i\nu)]_{(\nu=-i(n+1/2))}=\tilde{a}(n)=a_n$. We obtain:
\beq
\label{tredieci}
d_m=2\sqrt{\pi}\sum_{n=0}^\infty \frac{(-1)^n}{n!} a_n P_m\left [-i\left ( n+\frac{1}{2}\right )\right ].
\eeq
Next, we observe that:
\beq
\label{treundici}
\Gamma\left (\frac{1}{2}+i\nu\right )=\int_{-\infty}^{+\infty} e^{-i\nu v} e^{{-e}^{-v}} e^{-v/2} \, dv =
{\cF} \left\{ e^{{-e}^{-v}} e^{-v/2}\right\},
\eeq
where ${\cF}$ denotes the Fourier integral operator. Let us note that the function
$\exp(-\exp(-v)) \exp(-v/2)$ belongs to the Schwartz space $S$ of the $C^\infty(\R)$ functions which,
together with all their derivatives, decrease, for $|v|$ tending to $+\infty$,
faster than any negative power of $|v|$. Therefore we can write:
\beq
\label{tredodici}
\psi_m(\nu)=\frac{1}{\sqrt{\pi}}\,{\cF}\left\{P_m\left (-i\frac{d}{dv}\right )\left [e^{{-e}^{-v}} e^{-v/2}\right ]
\right\}.
\eeq
Substituting in expansion (\ref{treotto}) to the Pollaczek functions their representation
(\ref{tredodici}), we obtain:
\beq
\label{tretredici}
\tilde{a}\left (-\frac{1}{2}+i\nu\right )=\sum_{m=0}^\infty d_m\left\{\frac{1}{\sqrt{\pi}}\,{\cF}
\left [P_m\left (-i\frac{d}{dv}\right )\left [ e^{{-e}^{-v}} e^{-v/2}\right ]\right ]\right\}.
\eeq
Let us now apply the operator ${\cF}^{-1}$ to the r.h.s. of formula (\ref{tretredici}).
If we exchange the integral operator ${\cF}^{-1}$ with the sum, and this is legitimate within
the $L^2$-norm convergence, we obtain:
\beq
\label{trequattordici}
\begin{split}
& {\cF}^{-1} \sum_{m=0}^\infty d_m\left\{\frac{1}{\sqrt{\pi}}{\cF} \left [ P_m \left ( -i\frac{d}{dv}\right )
\left [ e^{{-e}^{-v}} e^{-v/2} \right ] \right  ] \right \} \\
&\qquad =\sum_{m=0}^\infty d_m\left\{\frac{1}{\sqrt{\pi}}{\cF}^{-1} {\cF} \left [ P_m\left (-i\frac{d}{dv}\right )
\left [e^{{-e}^{-v}} e^{-v/2}\right ] \right  ] \right \}.
\end{split}
\eeq
Finally, recalling formula (\ref{treuno}), we obtain the following expansion for the function
$\exp(v/2)F(v)$:
\beq
\label{trequindici}
e^{(v/2)}F(v)=\sum_{m=0}^\infty \frac{d_m}{\sqrt{\pi}} P_m\left (-i\frac{d}{dv}\right )
\left [e^{{-e}^{-v}} e^{-v/2}\right ],
\eeq
where the convergence is in the sense of the $L^2$--norm.

\noindent
Then it can be verified easily that:
\beq
\label{tresedici}
\sqrt{2}\, P_m\left (-i\frac{d}{dv}\right )\left\{e^{{-e}^{-v}} e^{-v/2}\right\} =
i^m\sqrt{2}\,L_m(2e^{-v})e^{{-e}^{-v}} e^{-v/2},
\eeq
where $L_m$ denotes the Laguerre polynomials (see also \cite{Koornwinder,Viano}).

It can be checked easily that the polynomials ${\cL}_m(v) \equiv i^m\sqrt{2}L_m(2e^{-v})$
are a set of polynomials orthonormal on the real line with the weight function
$w(v)=\exp(-v) \exp(-2\exp(-v))$, and, consequently, the set of functions
$\Phi_m(v)$, defined by formula (\ref{trequattro}) forms an orthonormal basis
in $L^2(-\infty,+\infty)$.
Finally, from formula (\ref{trequindici}) we obtain:
\beq
\label{trediciassette}
e^{v/2}F(v)= \sum_{m=0}^\infty c_m\left\{i^m\sqrt{2}L_m(2e^{-v})e^{{-e}^{-v}} e^{-v/2}\right\}=
\sum_{m=0}^\infty c_m\Phi_m(v) \qquad (v\in\R^+),
\eeq
where $c_m=d_m/\sqrt{2\pi}$, and the functions $\Phi_m(v)$ are given by
formula (\ref{trequattro}).
\end{proof}

The results of Theorem \ref{the:2} can be easily extended to all the cases considered in Section \ref{sezione2}.
For the sake of simplicity we limit ourselves to consider the Taylor series (\ref{dueventi})
treated in Theorem \ref{the:1prime}.

\begin{theorem}
\label{the:2prime}
If in the Taylor series $(\ref{dueventi})$ the coefficients $a_n$
satisfy condition $(\ref{unotre})$, then the jump function $F(x)$ can be represented
by the following series, which converges in the sense of the $L^2$--norm:
\beq
\label{ottantotto}
F(x)=\sum_{m=0}^\infty c_m \phi_m(x) \qquad (x\in (1,+\infty)).
\eeq
In the series $(\ref{ottantotto})$ the coefficients $c_m$ are given by formula
$(\ref{tretre})$, while the functions $\phi_m(x)$ are given by:
\beq
\label{ottantanove}
\phi_m(x) = i^m \sqrt{2} L_m\left (\frac{2}{x}\right ) \frac{e^{-1/x}}{x},
\eeq
$L_m$ being the Laguerre polynomials.
\end{theorem}

\begin{proof}
The proof of these results proceeds
exactly as in the case of the previous theorem. It is, indeed,
sufficient to observe that the Mellin transform (\ref{dueventidue}) can be easily
transformed in a Fourier--Laplace transform by putting: $x=e^v$. Let us remind that
we still denote (with a small abuse of language) the jump function by $F(x)$, as it
has already been noted before the formulation of Theorem \ref{the:1prime}. Accordingly, the basis
will be given by formula (\ref{ottantanove}). Finally, let us
note that the functions $\phi_m(x)$ form an orthonormal basis in $L^2(0,+\infty)$.
\end{proof}

Let us note that expansion (\ref{ottantotto}) furnishes a solution of the Cauchy
integral equation of the type (\ref{zeroquattroprimo}). However, let us observe that
up to now we have supposed that the coefficients $a_n$ are infinite in number and noiseless.
But this, in practice, is not the case. We have at our disposal only a finite number of
coefficients and, in addition, they are affected by noise or by round--off errors.
We are, therefore, forced to consider the following question: How to manage numerically
expansions (\ref{trediciassette}) and (\ref{ottantotto}).
Furthermore, let us note that the problem of reconstructing the jump function from
the coefficients $a_n$ is a classical example of ill--posed problem in the
sense of Hadamard \cite{Hadamard}. It is, indeed, strictly connected to the problem of the analytic
continuation up to the boundary of the analyticity domain.
In fact, let us focus our attention on the cut $z$--plane geometry considered in
Theorems \ref{the:1prime} and \ref{the:2prime}. In view of the Riemann mapping theorem, this cut plane can be
conformally mapped onto the unit disk in the $\zeta$--plane geometry (i.e., $|\zeta| < 1$)
through a suitable transformation $\zeta=\zeta(z)$. In this map the upper (lower) lip
of the cut is mapped in the upper (lower) half of the unit circle. Therefore, the
problem of solving the Cauchy--type integral equations (\ref{zeroquattroprimo})
corresponds to the analytic continuation up to the unit circle ($|\zeta|=1$). It is, then,
easy to exhibit Hadamard--like examples showing that the solution does not
depend continuously on the data in various topologies, including uniform and
$L^2$--topologies.

We shall treat all these questions in a separate paper devoted to the numerical
analysis. The main result that will be proved there reads as follows: If we take as
data a finite number $(N+1)$ of coefficients perturbed by noise $a_n^{(\epsilon)}$
(where $|a_n^{(\epsilon)} - a_n|\leq\epsilon,~n=0,1,2,\ldots,N,~\epsilon >0$)
we can still determine an approximation $F^{(\epsilon,N)}(x)$ of the jump function
$F(x)$, that asymptotically converges to $F(x)$, in the sense of the $L^2$--norm, as
$\epsilon\rightarrow 0$ and $N\rightarrow\infty$.

~

~

\subsection*{Acknowledgments}
One of us (G.A.V.) is deeply indebted to Prof. J. Bros for several illuminating discussions.

\end{document}